\documentclass{tran-l}

\begin{document}
\begin{center}
{\bf Heuristic algorithm for solving the graph isomorphism problem}\\
R.T. Faizullin, A.V. Prolubnikov\\
\end{center}
\ \\
\begin{quote} {\small In the paper we consider heuristic
algorithm for solving graph isomorphism problem. The algorithm
based on a successive splitting of the eigenvalues of the matrices
which are modifications (to positive defined) of graphs' adjacency
matrices. Modification of the algorithm allows to find a solution
for Frobenius problem. Formulation of the Frobenius problem is
following one. Given a pair of two matrices with the same number
of rows and columns. We must find out whether one of the matrix
can be acquired from another by permutation of it's rows and
strings or not. For example, solution of Frobenius problem can
give to us efficient way for decrypting of double permutation
cyphers problem for high dimension matrices.\\}
\end{quote}

\begin{center}
{\bf The graph isomorphism problem}
\end{center}
\ \\
\par Graph isomorphism problem is one of the problems for which we
can't say definitely whether this problem is polynomial or not
[1]. But it's known that this problem is polynomial for some
classes of graphs such as plane graphs, regular graphs and some
others [2], [3], [4]. Our algorithm is heuristic algorithm for
solving of the problem.

\par We consider heuristic algorithm for solving the graph
isomorphism problem. The algorithm based on successive splitting
of eigenvalues of the matrices that are modifications (to positive
defined) of graphs' adjacency matrices. During the work of the
algorithm we solving the linear equations that defines inverse
matrices of these matrices. Solutions of these systems gives to us
a permutation which is sought bijection.

\par There are two nonoriented graphs in the graph isomorphism
problem: graph $G_A=\langle V_A, E_A\rangle$ and graph $G_B =
\langle V_B, E_B\rangle$. $V_A,\ V_B$ are sets of graphs'
vertices. $E_A,\ E_B$ are sets of it's edges. They are supposed to
be the sets of the same power: $|V_A|=|V_B|,\ |E_A|=|E_B|$. Graph
isomorphism problem formulation is following one: whether exists
such bijection $\varphi : V_A\rightarrow V_B$ that if $(i,j)\in
E_A$ then $(\varphi(i),\varphi(j))\in E_B$ or not?

\par Algorithm works with modifications of adjacency matrices. Let $A_0$ be
an adjacency matrix of $G_A$, that is to say that $A_0=(a_{ij}^0)$
and
$$a_{ij}^0=
  \begin{cases}
    1,& \text{if}\ (i,j)\in E_A, \\
    0, & \text{else}.
  \end{cases}
$$
Let $B_0$ be an adjacency matrix of graph $G_B$.
\par Construct matrix $D_{A_0}$ according to $A_0$:
$$
\begin{pmatrix}
  d_1    &0      &\ldots & 0 \\
  0      &d_2    &\ldots & 0 \\
  \vdots &\vdots &\ddots & \vdots \\
  0      &0      &\ldots & d_n
\end{pmatrix}$$

$D_{A_0}$ is a scalar matrix with following elements:
$$d_i=\sum\limits_{j=1}^n a_{ij}^0+1.$$
Construct matrix $D_{B_0}$ similarly according to $B_0$. Let
$$A=A_0+D_{A_0},\ B=B_0+D_{B_0}. \eqno(1)$$ Matrices $A$ and $B$ are further
considerated matrices. They are matrices which algorithm works
with. They are symmetric and positive defined.

\par If graps $G_A=\langle V_A,\ E_A\rangle$ and $G_B=\langle V_B, E_B\rangle$
are isomorphic then we can acquire matrix $A$ from matrix $B$ by
successive permutations of it's rows with simultaneous
permutations of it's columns with same numbers. So we can
formulate the graph isomorphism problem as a particular case of
the Frobenius problem which formulation is following one. Can we
acquire matrix $B$ from matrix $A$ by successive permutations of
it's rows and columns?

\par Permutation of rows numbers $i$ and number $j$ for arbitrary matrix
is equal to right multiplying matrix by permutation matrix
$P_{ij}$. Permutation of columns is equal to left multiplying
matrix by the same permutation matrix. Permutation of vector's
components with numbers $i$ and $j$ take place for both left and
right multiplying vector by permutation matrix $P_{ij}$.

\par Let us consider the following systems of linear equations:
$$ Ax=e_j,\\ By=e_k. \eqno(2)$$ The vectors $e_i=(0,\ldots,0,1,0,\ldots,0)$
are basis vectors in the space $R^n$, matrices $A$ and $B$ such as
described above. Both systems are solvable and each solution is
unique solution because of $A$ and $B$ are matrices with diagonal
predominance and their determinants aren't equal to zero. Let
$x^j$ be a solution of the system $Ax=e_j$ and $y^k$ be a solution
of the system $By=e_k$.
\par Note that solution of the systems of linear equations (1) gives to
us inverse matrices of matrix $A$ and matrix $B$. So for $i$-th
component of vector $x^j$ the following is true:
$x_i^j=(-1)^{i+j}A_{ij}/\det(A)$, where $A_{ij}$ are algebraic
adjuncts for element $a_{ij}$ of the matrix $A$. That is to say
$x^j,\ j=1,\ldots ,n$ are columns of inverse matrix for $A$.

\par If $B=P_{jk}AP_{jk}$ then the following is true for solutions
 of the systems of linear equations (2):

 $$ x_i=y_i,\ i\neq j,\ i\neq k;$$
 $$x_j=y_k,\ x_k=y_j.$$
\par Indeed: $(Ax=e_j) \sim (P_{jk}Ax=P_{jk}e_j) \sim
(P_{jk}AxP_{jk}=P_{jk}e_jP_{jk}) \sim (P_{jk}AP_{jk}x=e_j) \sim
 (P_{jk}AP_{jk}xP_{jk}=e_jP_{jk}) \sim (BxP_{jk}=e_k).$ That is to say $xP_{jk}=y$.

\par If matrix $B$ can be acquired from matrix $A$ by successive permutations of it's rows with simultaneous
permutations of it's columns with same numbers then:
 $$B=P_{j_lk_l}\ldots P_{j_1k_1}AP_{j_1k_1}\ldots P_{j_lk_l},$$
and consequently $xP_{j_1k_1}\ldots P_{j_lk_l}=y$.

\par So if we'll change vector $e_k$ in the system of linear
 equations (2) with fixed $j$ by changing index $k$ from 1 to $n$ then
 vectors $x^j$ and $y^k$, which are corresponding solutions of the systems (2),
 will be the same vectors correct to permutation of their components only if
row with number $j$ of matrix $A$ corresponds to row with number
$k$ of matrix $B$. That is to say elements of row number $k$ of
matrix $B$ are permutated elements of row number $i$ of matrix
$A$. The same is true for columns of the matrices.

\par We'll accomplish the successive perturbations of its diagonal
 elements for finding unique row and unique column of the matrix $B$ that is
corresponds to row and column with number $i$ of the matrix $A$.
The algorithm works with symmetric matrices which can  be
transformed into scalar matrices by orthogonal transformations. So
$$\widetilde{A}=U_AAU_A^T,$$
$$ \widetilde{B}=U_BBU_A^T,$$
where $\widetilde{A},\ \widetilde{B}$ are scalar matrices with
eigenvalues on their diagonals and $U_A$, $U_B$ are matrices of
the orthogonal transformations. Diagonal elements of
$\widetilde{A}$ and $\widetilde{B}$ are eigenvalues of the
matrices $A$ and $B$.

\par The spectrums of the isomorphic graphs' adjacency matrices
are the same [5]. The same true for the matrices which algorithm
works with because of described above procedure of changing of the
diagonal elements of adjacency matrices leads only to shifting
spectrums of the matrices. So if the matrices are corresponds to
the isomorphic graphs then spectrums of the matrices will be
congruent after shifting.

\par It's obvious that if the spectrums of the both matrices is
simple or multiplicity of eigenvalues is not big then the problem
of graph isomorphism can be solved by a comparison of rows of the
matrices $U_A$, $U_B$, $\widetilde{A}$ and $\widetilde{B}$ [7].

\par The main difficulties arise with consideration of graphs which
spectrums contains multiple eigenvalues. Perturbation of the
matrices will allow us to perturb their spectrums while algorithm
works. Splitting of the eigenvalues will take place by this
perturbation. So we can maintain an unique correspondence between
rows and columns of the matrices.

\par If there is multiple eigenvalues in spectrums of
the matrices $À$ and $Â$ then they will be split and we'll be able
to find the permutation which is sought bijection $\varphi$. This
bijection establishes isomorphism of graphs $G_A = \langle V_A,\
E_A\rangle$ and $G_B = \langle V_B,\ E_B\rangle$. Computing
experiments give us evidence that split needed for determination
of the corresponding rows and columns take place much early than
at the last iteration of the algorithm. So at iteration with
number $\sqrt{n}$ there is no need for further perturbations of
diagonal elements of the matrices $A$ and $B$ for maintaining a
unique correspondence between vertices of graphs that represents
lattice on torus which numbers of vertices is equal to the $n$
where $n$ changes from 9 up to 400 .

\par At the algorithm's implementation on every iteration algorithm
works with already perturb matrices not with initial matrices. So
if we've got a correspondence between row number $j$ of matrix
$A^j$ and row number $k$ of matrix $B^j$ at the iteration number
$j$ and for columns with the same numbers then we'll considerate
further perturb matrices $A^{j+1}$ and $B^{j+1}$:
$$A^{j+1}=A^j+\varepsilon C^j,\  B^{j+1}=B^j+\varepsilon C^k.$$

We make perturbation by scalar matrices $C^k$ with elements $c_i$
on it's diagonal:
$$c_i=
  \begin{cases}
    1, & \text{if}\ i=j=k, \\
    0, & \text{else}.
  \end{cases}$$

\par As a result if matrix $A$ can be acquired from matrix $B$ by successive
permutations of it's rows with simultaneous permutations of it's
columns with same numbers then we get the sought permutation $P$
while $j$ changes from $1$ to $n$.  That is to say
$$P=\begin{pmatrix}
  1   & 2   & \ldots & n \\
  k_1 & k_2 & \ldots & k_n
\end{pmatrix},$$
$k_j$ is a number of row of matrix $B$ that obtained at $j$-th
iteration of the algorithm.

\par As a matter of fact permutation $P$ is one of the possible permutations.
So $P$ is a bijection $\varphi : V_A \rightarrow V_B$ that sets
isomorphism of graphs $G_A$ and $G_B$.

\smallskip

\begin{center}
{\bf The spectral splitting algorithm}\\
(First schema)
\end{center}
Step 0. $A^0:=A$,\ $j:=1$.\\
Step 1. If $j<n$ then
$A^{j}:=A^{j-1}+\varepsilon C^j$, else stop the algorithm's implementation.\\
Step 2. Solving of the system of linear equations $Ax=e_j$. $x^j$ is the solution.\\
Step 3. $k:=1$. If $k<n$ then go to Step 3.1, else go to Step 4.\\
Step 3.1. $B^k:=B^{k-1}+\varepsilon C^k$.\\
Step 3.2. Solving of the system of linear equations $B^ky=e_k$. $y^k$ is the solution\\
Step 3.3. $k:=k+1$. Go to Step 3.\\
Step 4. Comparing norms of $x^j$ and $y^k$, $k=1,\ldots ,n$.\\
If $\forall k\ ||x^j||\neq||y^k||$ then graphs $G_A$ and $G_B$
are not isomorphic. Stop the algorithm's implementation.\\
If exists $k: ||x^j||=||y^k||$ è $x_j^j=y_k^k$ then $P(j):=k$
(Maintaining a correspondence between vertex $j$ of graph $G_A$
and vertex $k$ of graph $G_B$), $B^j:=B^{j-1}+\varepsilon C^k$.
Else graphs $G_A$ and $G_B$ are not isomorphic. Stop the algorithm implementation.\\
Step 5. $j:=j+1$. Go to Step 1.

\bigskip
\par The hardness of this scheme is equal to $O(n^4)$,
$n$ is a number of rows at the square matrices $A$ and $B$. We
have to notice that this scheme can be modified to another scheme
which laboriousness is equal to $O(n^{3.5})$. We can make it by
adding the split checking procedure (Step 6 and step 7 of the
scheme stated below).

\begin{center}
{\bf The spectral splitting algorithm}\\
(Second schema)
\end{center}
Step 0. $A^0:=A$,\ $j:=1$.\\
Step 1. If $j<n$ then
$A^{j}:=A^{j-1}+\varepsilon C^j$, else stop the algorithm's implementation.\\
Step 2. Solving of the system of the linear equations $Ax=e_j$. $x^j$ is the solution.\\
Step 3. $k:=1$. If $k<n$, then go to Step 3.1, else then go to Step 4.\\
Step 3.1. $B^k:=B^{j-1}+\varepsilon C^k$.\\
Step 3.2. Solving of the system of the linear equations $B^ky=e_k$. $y^k$ is the solution.\\
Step 3.3. $k:=k+1$. Go to Step 3.\\
Step 4. Comparing norms of $x^j$ and $y^k$, $k=1,\ldots ,n$.\\
If $\forall k\ ||x^j||\neq||y^k||$ then graphs $G_A$ and $G_B$ are not isomorphic. Stop the algorithm's implementation.\\
If exists $k: ||x^j||=||y^k||$ è $x_j^j=y_k^k$ then $P(j):=k$
(Maintaining a correspondence between vertex $j$ of graph $G_A$
and vertex $k$ of graph $G_B$), $B^j:=B^{j-1}+\varepsilon C^k$.
 Else graphs $G_A$ and $G_B$ are not isomorphic. Stop the algorithm's implementation.\\
Step 5. $j:=j+1$.\\
Step 6. Solving of the system of linear equations $A^jx=e_l$. $x^l,\ l=1,\ldots , n$ is the solution.\\
Step 7. If $\forall l\ \forall p: ||x^l||\neq||x^p||$ Go to Step 8, else go to Step 1.\\
Step 8. If $j<n$ then implement solving of the system of the linear equations $B^ky=e_k,\ k=1,\ldots ,n;$ where $k$ such that not exists $ i: P (i)=j$.\\
Step 9. Comparing norms of $||x^l||$ and $||y^k||,\ l=j,\ldots
,n,\ k=1,\ldots ,n;$ and $k$ such that not exists $i: P (i)=j$. If
$||x^l||=||y^k||$ then $P(l):=k$.\\
Step 10. Stop the algorithm's implementation.

\par Though the spectral splitting algorithm show oneself good by
computing experiments we have to notice that we can't guarantee
finding solution of the graph isomorphism problem for arbitrary
graphs because of we can't say that the following situation is
impossible. There are two vectors $x_j$ and $y_k$ with equal or
such close norms for we can't say whether one vector can be
acquired from another by permutation of it's components or not
because of we implements check procedure only by one components of
each vectors. They are components $x_j^j$ and $y_k^k$.
\par In particular we realized computing experiments for regular graphs with number of vertices
up to 2500 for the graph isomorphism problem.

\bigskip

\begin{center}
{\bf Solving of the Frobenius problem for arbitrary\\ square
matrices of complete rank}
\end{center}

\par The spectral splitting algorithm for the graph isomorphism
problem can be applied for weighted graphs without any
modifications. Scheme of the algorithm is the same; the only
difference in matrices algorithm works with.

\par There are two nonoriented graphs in the graph isomorphism
problem: graph $G_A=\langle V_A, E_A\rangle$ and graph $G_B =
\langle V_B, E_B\rangle$. $V_A,\ V_B$ are sets of graphs'
vertices, $E_A,\ E_B$ are sets of it's edges. They are supposed to
be the sets of the same power: $|V_A|=|V_B|,\ |E_A|=|E_B|$. Given
a function $H_A:E_A\rightarrow R$ è $H_B:E_B\rightarrow R$ which
defines weights of graphs' edges. Graph isomorphism problem for
weighted graphs formulation is following one: whether exists such
bijection $\varphi : V_A\rightarrow V_B$, that if $(i,j)\in E_A$,
then $(\varphi(i),\varphi(j))\in E_B$ and if $(i,j)\in E_A$ then
$(\varphi(i),\varphi(j))\in E_B$ è $H_A(i,j)=H_B(\varphi
(i),\varphi (j))$ or not?

\par The adjacency matrices of weighted nonoriented graphs
transforms into positive defined matrices with diagonal
predominance, but we should select its diagonal elements in such
way that conditionality number should be confined.

\par Let matrix $A_0=(a_{ij}^0)$ be an adjacency matrix of
weighted nonoriented graph $G_A$ and matrix $D_{A_0}$ be a
diagonal matrix with such elements $d_i$ that:
$$d_i=d+\sum\limits_{j=1}^na_{ij}^0,$$
where $d=\max\limits_{1\leq i\leq n}\sum\limits_{j=1}^na_{ij}^0$.

Construct matrix $D_{B_0}$ for weighted nonoriented graph $G_B$ by
same way according to adjacency matrix of graph $B_0$.

\par For conditionality number $\mu (A)$ of the symmetric matrix $A$
the following is true [6]:
$$\mu (A)\leq \displaystyle\frac {\eta (A)}{\chi (A)},$$
where $\eta (A)=\max\limits_{1\leq i\leq
n}(a_{ii}+\sum\limits_{j\neq i}|a_{ij}|),\ \chi
(A)=\min\limits_{1\leq i\leq n}(a_{ii}-\sum\limits_{j\neq
i}|a_{ij}|)$.
\par By our choice of $d$:\\
$\eta (A)=\max\limits_{1\leq i\leq
n}(a_{ii}+\sum\limits_{j=1}^n|a_{ij}^0|)=
a_{i_1i_1}+\sum\limits_{j=1}^n|a_{i_1j}^0|=d_{i_1}+\sum\limits_{j=1}^n|a_{i_1j}^0|=d+2\sum\limits_{j=1}^n|a_{i_1j}^0|=$\\
$=3\sum\limits_{j=1}^n|a_{i_1j}^0|$,\\
$\chi (A)=\min\limits_{1\leq i\leq
n}(a_{ii}-\sum\limits_{j=1}^n|a_{ij}^0|)=a_{i_2i_2}-\sum\limits_{j=1}^n|a_{i_2j}^0|=
d_{i_2}-\sum\limits_{j=1}^n|a_{i_2j}^0|=\sum\limits_{j=1}^n|a_{i_2j}^0|+d-$\\
$-\sum\limits_{j=1}^n|a_{i_2j}^0|=\sum\limits_{j=1}^n|a_{i_2j}^0|+\sum\limits_{j=1}^n|a_{i_1j}^0|
-\sum\limits_{j=1}^n|a_{i_2j}^0|=\sum\limits_{j=1}^n|a_{i_1j}^0|.$\\
Consequently
$$
\mu (A)=\displaystyle\frac{\eta (A)}{\chi
(A)}\leq\displaystyle\frac{3\sum\limits_{j=1}^n|a_{i_1j}^0|}
{\sum\limits_{j=1}^n|a_{i_1j}^0|}=3.$$

\par There are two matrices $F_A$ and $F_B$ with the same number of rows and columns in the Frobenius
problem. We must find out whether one matrix can be acquired from
another by some permutation of it rows and columns.

\par Suppose that matrices $F_A$ and $F_B$ are square matrices of the complete rank
with number of rows is equal to $n$. Matrix $A_0$ of the following
structure:
$$\begin{pmatrix}
  0     & F_A \\
  F_A^T & 0
\end{pmatrix}.$$
This matrix can be considered as an adjacency matrix of some
weighted nonoriented graph. Construct matrix $D_{A_0}$ with
structure described above for matrix $A_0$. Then construct matrix
$A$:
$$A=A_0+D_{A_0}.$$
That is to say
$$A=\begin{pmatrix}
  D_{A_0}^1 & F_A \\
  F_A^T     & D_{A_0}^2
\end{pmatrix},$$
and
$$D_{A_0}=\begin{pmatrix}
  D_{A_0}^1 & 0 \\
  0       & D_{A_0}^2
\end{pmatrix},$$
where
$$D_{A_0}^1=\begin{pmatrix}
d_1    & 0      & \ldots &  0 \\
0      & d_2    & \ldots &  0 \\
\vdots & \vdots & \ddots &  \vdots \\
0      & 0      & \ldots &  d_n\\
\end{pmatrix},\
D_{A_0}^2=\begin{pmatrix}
d_{n+1}    & 0      & \ldots &  0 \\
0      & d_{n+2}    & \ldots &  0 \\
\vdots & \vdots & \ddots &  \vdots \\
0      & 0      & \ldots &  d_{2n}\\
\end{pmatrix}.
$$
So matrix $A$ is matrix for initial matrix $F_A$ for algorithm to
work with.
\par Matrix $A$ is positive defined symmetric matrix with diagonal
predominance. It corresponds to some weighted nonoriented
dicotyledonous graph. Matrix $F_B$ is a square matrix which number
of rows the same as the number of rows of matrix $F_A$.

Construct matrix $B$ by the same way for matrix $F_B$. Form the
solving permutation of rows and columns of matrices $A$ and $B$ by
applying the spectral splitting algorithm to matrices $A$ and $B$.

\par Form the solving permutation for initial matrices $F_A$ and
$F_B$ by the permutations of rows and columns of the matrices $A$
and $B$ in the following way. Let it given by the algorithm that
the row and column with number $i$ of matrix $A$ corresponds to
the row and column with number $j$ of matrix $B$. It is obvious
that if $1\leq i\leq n,\ 1\leq j\leq n$ then the permutation of
pair of rows $\{i,j\}$ with simultaneous permutation of it's
columns $\{i,j\}$ of matrix $A$ corresponds to permutation of rows
number $i$ and number $j$ of initial matrix $F_A$. If $n+1\leq
i\leq 2n,\ n+1\leq j\leq 2n$ then permutation of pair of rows
$\{i,j\}$ with simultaneous permutation of columns $\{i,j\}$ of
matrix $A$ corresponds to permutation of columns with number $i$
and number $j$ of initial matrix $F_A$. As soon as matrices $A$
and $B$ which algorithm works with have the structure as following
one:
$$\begin{pmatrix}
  *      & 0      & \ldots & 0      & *      & *      & \ldots & *  \\
  0      & *      & \ldots & 0      & *      & *      & \ldots & *  \\
  \vdots & \vdots & \ddots & \vdots & \vdots & \vdots &        & \vdots \\
  0      & 0      & \ldots & *      & *      & *      & \ldots & *  \\
  *      & *      & \ldots & *      & *      & 0      & \ldots & 0 \\
  *      & *      & \ldots & *      & 0      & *      & \ldots & 0 \\
  \vdots & \vdots &        & *      & \vdots & \vdots & \ddots & \vdots \\
  *      & *      & \ldots & *      & 0      & 0      & \ldots & *
\end{pmatrix},$$
($*$ are positions admissible for nonzero elements) so the
situation of maintaining of the correspondence for pair of rows
and pair of columns of matrices $A$ and $B$ with such numbers $i$
and $j$ that $\{i,\ j\}$ and $1\leq i\leq n,\ n+1\leq j\leq 2n$ is
impossible. It impossible because of it means appearance of
nonzero nondiagonal elements at the left and upper part of the
matrix $A\ (\text{of the matrix}\ D_{A_0}^1)$ and right and lower
part $(\text{of the matrix}\ D_{A_0}^2)$ if only the initial
matrices $F_A$ and $F_B$ are not scalar. This is impossible by
constructing of the matrices $A$ and $B$. Consequently appearance
of such correspondence is impossible while algorithm works.

\bigskip
\begin{center}
{\bf Deciphering of the double permutation cipher\\ by the
spectral splitting algorithm}
\end{center}
\par Suppose we have a text presented by some square matrix of symbols.
Let us form matrix of the cipher with some digital code for every
symbol of text. The second matrix presented some coding by this
way text too. It's needed to find out whether the text presented
by the second matrix is the double permutation cipher for the firs
matrix or not. That is to say whether this text can be obtained
from the first text presented by matrix by permutation of it's
rows and columns. Such problem may be interpreted as a case of the
Frobenius problem for digital matrices. So we can apply the
spectral splitting algorithm for solving of this problem.

\par There are computational experiments was implemented for deciphering of the double permutation
ciphers for texts with number of symbols up to $10000$.

\begin{center}
{\bf References}
\end{center}

1. M. Garey and D. Jhonson {\it Computers and Intractability, A
Guide to Theory of NP-completness.} Freeman and Co., 1979\\

2. Hopcroft, Wong {\it A linear time algorithm for isomorphism of
planar graphs.} Proceedings of the Sixth Annual ACM Symposium on
Theory of Computing, p. 172-184, 1974.\\

3. Luks {\it Isomorphism of graphs of bounded valence can be
tested
in polynomial time.} Proc. 21st IEEE FOCS Symp., 42,\ 49, 1980.\\

4. Hoffmann  {\it Group-Theoretic Algorithms and Graph
Isomorphism.}
Lecture Notes in Computer Science (Chapter V). P.127-138, 1982\\

5. Tzvetkovich D. and others {\it Spectrums of graphs. Theory and
application.} Kiev: The scientific thought, 1984.\\

6. Godunov S. and others  {\it Guaranteed exactness of solutions
of systems of linear equations in euclidian spaces.} Novosibirsk:
Science, 1988.\\

7. Kikina A., Faizullin R. {\it The algorithm for testing graph
isomorphism.} VINITI 21.06.95 1789-Â95.\\

\end{document}